\input amstex
\magnification=\magstep1 
\baselineskip=13pt
\documentstyle{amsppt}
\vsize=8.7truein \CenteredTagsOnSplits \NoRunningHeads
 \def\EE{\bold{E}\thinspace}
 \def\PP{\bold{P}\thinspace}
 
 \def\tr{\operatorname{trace}}
 \def\per{\operatorname{per}}
 \def\ham{\operatorname{ham}}
 \def\wei{\operatorname{weight}}
 \def\spt{\operatorname{spt}}
 \topmatter
 
\title Partition functions for dense instances of combinatorial enumeration problems \endtitle
\author Alexander Barvinok \endauthor
\address Department of Mathematics, University of Michigan, Ann Arbor,
MI 48109-1043, USA \endaddress
\email barvinok$\@$umich.edu \endemail
\date May 2013 \enddate
\thanks  This research was partially supported by NSF Grant DMS 0856640.
\endthanks 
\abstract Given a complete graph with positive weights on its edges, we define the weight of a subset of edges as the product of weights of the edges in the subset and consider sums (partition functions) of weights over subsets of various kinds: cycle covers, closed walks, spanning trees. We show that if the weights of the edges of the graph are within a constant factor, fixed in advance, of each other then the bulk of  the partition function is concentrated on the subsets  of a particularly simple structure: cycle covers with few cycles, walks that visit every vertex only few times, and spanning trees with small degree of every vertex. This allows us to construct a polynomial time algorithm to separate graphs with many Hamiltonian cycles from graphs that are sufficiently far from Hamiltonian.
\endabstract
\endtopmatter

\document

\head 1. Introduction and main results \endhead

Given a graph $G=(V, E)$ with set $V$ of vertices and set $E$ of edges, it is a classical NP-complete problem to determine whether $G$ is Hamiltonian. A Hamiltonian cycle in $G$ may be alternatively described as a cycle cover consisting of a single cycle or a closed walk which visits every vertex once, while a Hamiltonian path in $G$ can be described as a spanning tree with the degree of every vertex not exceeding 2. Let us embed $G$ into the complete graph $K=(V, E_1)$ with the same set $V$ of vertices and assign weights to the edges of $K$ by 
$$w(e)=\cases 1 &\text{if\ } e \in E \\ \delta &\text{if \ } e \notin E, \endcases$$
where $0 < \delta < 1$ is a small positive number. For a subset $S \subset E_1$, we define the weight of $S$ by 
$$w(S) = \prod_{e \in E_1} w(e)$$
and consider the sum ({\it partition function})
$$f({\Cal C}, w)=\sum_{S \in {\Cal C}} w(S),$$
where ${\Cal C}$ is the class of sets of interest: cycle covers in $K$, closed walks with $|V|$ steps in $K$ or spanning trees in 
$K$. Thus $w(S)=1$ if and only if $S$ consists of edges of the original graph $G$, and the more non-edges of $G$
the set $S$ uses, the smaller its contribution $w(S)$ towards $f({\Cal C}, w)$. This paper is based on the following two observations: 
\medskip
1) it is easy to compute $f({\Cal C}, w)$ when ${\Cal C}$ is the class of cycle covers, or closed walks of a given length, or spanning trees in $K$ 
\smallskip
and
\smallskip
 2) if $\delta >0$ is fixed in advance or does not decrease too fast with the size of the problem, the bulk of $f({\Cal C}, w)$ is 
contributed by sets $S$ of a particularly simple structure: cycle covers with $O\bigl(\ln |V|\bigr)$ cycles, closed walks with $|V|$ steps that do to visit any vertex more than $O\bigl(\ln |V|/\ln \ln |V|\bigr)$ times and spanning trees where the degree of every vertex is $O\bigl(\ln |V|/\ln \ln |V|\bigr)$. 
\medskip
This allows us to use easily computable partition functions (such as those corresponding to cycle covers) to approximate partition functions that are hard to compute  (such as those corresponding to Hamiltonian cycles). In particular, we obtain a polynomial time algorithm to separate graphs that have many Hamiltonian cycles (at least $\epsilon^{|V|} |V|!$ for some fixed $0 < \epsilon < 1$) from graphs that are far from Hamiltonian (where each Hamiltonian cycle in $K$ contains at least $\gamma |V|$ non-edges of $G$, for some fixed $0< \gamma < 1$).

\subhead (1.1) Permanents and Hamiltonian permanents \endsubhead
Let $A=\left(a_{ij}\right)$ be an $n \times n$ real matrix. The {\it permanent of $A$} is defined as
$$\per A =\sum_{\sigma \in S_n} \prod_{i=1}^n a_{i \sigma(i)},$$
where the sum is taken over the symmetric group $S_n$ of permutations of the set $\{1, \ldots, n\}$. As is known, the problem
of computing the permanent exactly is $\#P$-hard, even if the entries of $A$ are restricted to be 0 and 1 \cite{Va79}. For non-negative matrices a fully polynomial randomized approximation scheme is available \cite{J+04}. We, however, are interested in computing permanents of a rather restricted class of matrices. Namely, let us fix a $\delta >0$ and suppose that 
$$\delta \ \leq \ a_{ij} \ \leq \ 1 \quad \text{for all} \quad i, j. \tag1.1.1$$
Then the scaling algorithm of \cite{L+00}, see also \cite{BS11}, approximates $\per A$ within a factor of $n^{O(1)}$, where the implied constant in the $O(1)$ notation depends on $\delta$. The advantage of the algorithm of \cite{L+00} is that beside being polynomial time, it is deterministic and easy to implement.

Let $H_n \subset S_n$ be the subset of $(n-1)!$ permutations consisting of a single cycle. We define the {\it Hamiltonian permanent} by 
$$\ham A =\sum_{\sigma \in H_n} \prod_{i=1}^n a_{i \sigma(i)}.$$
If $A$ is a 0-1 matrix then it is an NP-complete problem to tell $\ham A$ from $0$, as the problem is equivalent to testing Hamiltonicity of the directed graph with the adjacency matrix $A$. It turns out, however, that when (1.1.1) holds, $\per A$ and $\ham A$ have the same logarithmic order.

\proclaim{(1.2) Theorem} Let us fix a $0 < \delta < 1$. Then there exists a $\gamma=\gamma(\delta)>0$ such that for any 
$n \times n$ matrix $A=\left(a_{ij}\right)$ which satisfies (1.1.1), we have 
$$(\delta n)^{-\gamma \ln n} \per A \ \leq \ \ham A \ \leq \ \per A.$$
\endproclaim

Our proof of Theorem 1.2 is based on the following observation: under the condition (1.1.1), the bulk of the terms in $\per A$ falls on the permutations $\sigma \in S_n$ with few cycles.

\proclaim{(1.3) Theorem} For $\sigma \in S_n$ let $c(\sigma)$ denote the number of cycles in $\sigma$. Let $A=\left(a_{ij}\right)$ be an $n \times n$ matrix such that (1.1.1) holds. Then 
$$\sum \Sb \sigma \in S_n: \\ c(\sigma) \leq 4 +4 \delta^{-2} \ln n \endSb \prod_{i=1}^n a_{i \sigma(i)} \ \geq \ 
{1 \over 2} \per A.$$
\endproclaim

In a somewhat different setting, the relation between the permanent and Hamiltonian permanent of the adjacency matrix of a $k$-regular graph was used in \cite{Vi12}, see also \cite{SV13}.

Theorem 1.3 describes what appears to be a fairly general phenomenon: the partition function on dense instances concentrates on objects of a simple structure. We give two more examples. 

\subhead (1.4) Walks in a graph \endsubhead 
A {\it closed walk} $\pi$ in the complete directed graph with vertices $1, \ldots, n$ is just a sequence 
$$\pi=i_1 \rightarrow i_2 \rightarrow \ldots \rightarrow i_n \rightarrow i_1 \tag1.4.1$$
of not necessarily distinct numbers $i_1, \ldots, i_n \in \{1, \ldots, n\}$. Given an $n \times n$ matrix $A=\left(a_{ij}\right)$, we define the {\it weight} of the walk (1.4.1) by 
$$\wei(\pi)=a_{i_1 i_2} a_{i_2 i_3} \cdots a_{i_{n-1} i_n} a_{i_n i_1}.$$
Let $\Pi_n$ be 
the set of all $n^n$ closed walks of length $n$.
Then
$$\tr A^n = \sum_{\pi \in \Pi_n} \wei(\pi).  \tag1.4.2$$
We define the {\it degree} of a vertex $i$ in a walk $\pi$ as the number $\deg_i(\pi)$ of times the walk arrives to $i$, that is, the number of steps $\ast \rightarrow i$ in (1.4.1). For example, a Hamiltonian cycle is a walk $\pi$ such that $\deg_i(\pi)=1$ for $i=1, \ldots, n$. It turns out that if $A$ satisfies (1.1.1) then the bulk of (1.4.2) falls on the walks with small (sublogarithmic) degrees of the vertices.
\proclaim{(1.5) Theorem} Let $A=\left(a_{ij}\right)$ be an $n \times n$ matrix which satisfies (1.1.1). 
Then 
$$\sum \Sb \pi \in \Pi_n: \\ \deg_i(\pi) \leq 3\ln n/\delta^2 \ln \ln n \\ \text{for} \ i=1, \ldots, n \endSb \wei(\pi) \ \geq \ 
\left({n-1 \over n}\right) \tr A^n.$$
\endproclaim

\subhead (1.6) Spanning trees in a graph \endsubhead Let us consider the complete undirected graph on $n$ vertices 
$\{1, \ldots, n\}$, without loops or multiple edges and let $T_n$ be the set of all $n^{n-2}$ {\it spanning trees} in the graph.
Given an $n \times n$ real symmetric matrix $A=\left(a_{ij}\right)$, we define the {\it weight} of a spanning tree $\tau$ as 
$$\wei(\tau)=\prod_{\{i,j\} \text{\ is an edge of\ } \tau} a_{ij}.$$
We define the partition function of spanning trees by 
$$\spt A = \sum_{\tau \in T_n} \wei(\tau). \tag1.6.1$$
As is well known, the Kirchoff formula gives a fast algorithm of computing $\spt A$. Namely, we orient the edges of the complete graph arbitrarily, consider the $ n \times {n \choose 2}$ incidence matrix $B$, with rows indexed by vertices $\{1, \ldots, n\}$, 
columns indexed by directed edges $e=(i \rightarrow j)$ and entries
$$b_{ie}=\cases a_{ij} &\text{if\ } e=i \rightarrow j \\ -a_{ij} &\text{if\ } e=j \rightarrow i \\ 0 &\text{otherwise.} \endcases$$
If $\widehat{B}$ is obtained from $B$ by crossing out an arbitrary row, then 
$$\spt A = \det B^T B,$$
see, for example, Section II.3 of \cite{Bo98}.
Again, it turns out that once (1.1.1) holds, the bulk of (1.6.1) falls on the trees with small degrees of vertices. Denoting $\deg_i(\tau)$ the degree of vertex $i$ in the spanning tree $\tau$, we obtain the following result.

\proclaim{(1.7) Theorem} Let $A=\left(a_{ij}\right)$ be an $n \times n$ matrix which satisfies (1.1.1). Assuming that 
$n \geq 2/(1-\delta)$, we have
$$\sum \Sb \tau \in T_n: \\ \deg_i (\tau) \leq  1+3\ln n/\delta \ln \ln n  \\ \text{for\ } i=1, \ldots, n \endSb \wei(\tau) \ \geq \ 
\left({n-1 \over n} \right) \spt A.$$
\endproclaim 

We obtain Theorems 1.5 and 1.7 as a corollary to the following general result.

\proclaim{(1.8) Theorem} For positive integers $m$ and $n$, let $\Delta_{m, n}$ be the set of non-negative integer 
vectors $\left(\alpha_1, \ldots, \alpha_n \right)$ such that $\alpha_1 + \ldots  + \alpha_n=m$. 
Suppose that there is a probability measure on $\Delta_{m, n}$ such that 
$$\PP\bigl(a\bigr) = {w\left(\alpha_1, \ldots, \alpha_n\right) \over \alpha_1! \cdots \alpha_n!} \quad \text{where} \quad 
a=\left(\alpha_1, \ldots, \alpha_n\right),$$
and for some $0 < \delta \leq 1$ non-negative numbers $w(a)$ satisfy 
$$ w\left(\alpha_1, \ldots, \alpha_n \right)  \ \leq \  \delta^{-1} w\left(\beta_1, \ldots, \beta_n \right) \quad \text{whenever} \quad
\sum_{i=1}^n \left| \alpha_i - \beta_i\right| =2. $$
If $m \geq \delta n$ then 
$$\PP\left(\left(\alpha_1, \ldots, \alpha_n\right): \ \max_{i=1, \ldots, n} \alpha_i \ \geq \ {3 m \ln n \over \delta n \ln \ln n} \right) 
\ \leq \ {1 \over n}.$$
\endproclaim

\subhead (1.9) Applications to testing Hamiltonicity of graphs \endsubhead
Let $G=(V, E)$ be a finite directed graph without loops or multiple edges directed the same way, with set $V$ of vertices and  set $E$ of edges.
We identify $V=\{1, \ldots, n\}$, after which $G$ is represented by its adjacency matrix
$A=\left(a_{ij}\right)$, where 
$$a_{ij}=\cases 1 &\text{if\ } (i \rightarrow j) \in E \\ 0 &\text{otherwise.}\endcases$$
Then $\ham A$ is the number of Hamiltonian cycles in $G$.
Let us fix constants $0 < \epsilon, \gamma < 1$. 
We want to distinguish the following two cases:
\medskip
a) the graph $G$ has at least $\epsilon^n (n-1)!$ Hamiltonian cycles 
\smallskip
\noindent and 
\smallskip
b) any Hamiltonian cycle in the complete graph with the set $V$ of vertices contains at least $\gamma n$ non-edges of $G$.
\medskip
Let us choose a positive $\delta < \epsilon^{1/\gamma}$ and consider a perturbation $B=\left(b_{ij}\right)$ of the matrix $A$, defined as follows:
$$b_{ij}=\cases 1 &\text{if\ } a_{ij} =1 \\ \delta &\text{if \ } a_{ij}=0. \endcases$$
If a) holds we have, obviously,
$$\ham B \ \geq \ \ham A \ \geq \ \epsilon^n (n-1)!. \tag1.9.1$$
On the other hand, if b) holds,  then 
$$\delta^{\gamma n} (n-1)! \ \geq \ \ham B. \tag1.9.2$$
Comparing (1.9.1) and (1.9.2) and using that by Theorem 1.2
$$\ham B = n^{O(\ln n)} \per B,$$
with the implicit constant in the ``$O$" notation depending on $\delta$, 
we conclude that we can distinguish in polynomial time between the alternatives a) and b) for any fixed $\epsilon$ and $\gamma$.

Similarly, one can separate in polynomial time graphs containing many Hamiltonian cycles from graphs that don't have closed walks or spanning trees with small degrees of vertices.

\head 2. Proofs of Theorems 1.2 and 1.3 \endhead

Let us fix a positive $n \times n$ matrix $A=\left(a_{ij}\right)$. We consider $S_n$ as a finite probability space, where we let
$$\PP(\sigma) =\left( \per A \right)^{-1} \left( \prod_{i=1}^n a_{i \sigma(i)} \right) \quad \text{for} \quad \sigma \in S_n.$$
\proclaim{(2.1) Lemma} Let us define random variables
$$l_i:\ S_n \longrightarrow {\Bbb R} \quad \text{for} \quad i=1, \ldots, n,$$
where $l_i(\sigma)$ is the length of the cycle of $\sigma$ containing $i$. Suppose that (1.1.1) holds.
Then
$$\PP\left(\sigma:\ l_i(\sigma)=m \right) \ \leq \ {1 \over \delta^2(n-m)} \quad \text{for} \quad i=1, \ldots, n$$
and $m=1, \ldots, n$.
\endproclaim 
\demo{Proof} Without loss of generality, we assume that $i=1$. With  the set of permutation $\sigma \in S_n$ such that 
$l_1(\sigma)=m$ we associate a set $\Sigma \subset S_n$ as follows. We write the cycle of $\sigma$ containing 1
as
$$1 =j_1 \rightarrow j_2 \rightarrow \ldots \rightarrow j_m \rightarrow 1.$$
Let us pick any of the $n-m$ numbers, say $r$, not in the cycle. We write the cycle containing $r$ as 
$$r=j_{m+1} \rightarrow j_{m+2} \rightarrow \ldots \rightarrow j_{m+k} \rightarrow r$$
and produce a permutation $\tau \in \Sigma$ by merging the two cycles together:
$$1=j_1 \rightarrow j_2 \rightarrow \ldots \rightarrow j_m \rightarrow r=j_{m+1} \rightarrow j_{m+2} \rightarrow \ldots \rightarrow j_{m+k} \rightarrow 1.$$ Because of (1.1.1), we have
$$\PP(\tau) \ \geq \ \delta^2 \PP(\sigma). \tag2.1.1$$
The set $\Sigma$ consists of all permutations $\tau$ thus obtained from all permutations $\sigma$ with $l_1(\sigma)=m$.
We observe that every $\tau \in \Sigma$ is obtained from a unique permutation $\sigma$. To reconstruct $\sigma$ from $\tau$, 
we choose the cycle of $\tau$ containing 1, write it as 
$$1 \rightarrow j_1 \rightarrow j_2 \rightarrow \cdots \rightarrow j_{m+k} \rightarrow 1$$
for some $k>0$ and split it into the two cycles,
$$1 \rightarrow j_1 \rightarrow j_2 \rightarrow \ldots \rightarrow j_m \rightarrow 1 \quad \text{and} \quad 
j_{m+1} \rightarrow \ldots \rightarrow j_{m+k} \rightarrow j_{m+1}.$$
Since every permutation $\sigma \in S_n$ with $l_1(\sigma)=m$ gives rise to $n-m$ permutations $\tau \in \Sigma$, using 
(2.1.1) we obtain
$$\PP\left(\sigma:\ l_1(\sigma)=m\right) \ \leq \ {1 \over \delta^2 (n-m)} \PP\left(\tau:\ \tau \in \Sigma \right) \ \leq \ 
{1 \over \delta^2(n-m)},$$
as desired.
{\hfill \hfill \hfill} \qed
\enddemo

\proclaim{(2.2) Lemma} Let us consider a random variable 
$$c\: S_n \longrightarrow {\Bbb R},$$
where $c(\sigma)$ is the number of cycles of a permutation $\sigma \in S_n$. Suppose that (1.1.1) holds.
Then 
$$\EE c \ \leq \ 2 + 2\delta^{-2} \ln n.$$
\endproclaim
\demo{Proof} Let $l_i$ be the random variables of Lemma 2.1. Then 
$$c=\sum_{i=1}^n l_i^{-1},$$
since for any $\sigma \in S_n$ the sum of $1/l_i(\sigma)$ for all $i$ in a cycle of $\sigma$ of length $l$ is 1.
Using Lemma 2.1, we obtain
$$\split \EE  l_i^{-1} = &\sum_{m=1}^n {1 \over m} \PP\left(\sigma: \ l_i(\sigma)=m \right) \\=&\sum_{m:\ m \leq n/2} {1 \over m}
\PP\left(\sigma:\ l_i(\sigma)=m\right) + \sum_{m:\  m> n/2} {1 \over m} \PP\left(\sigma:\ l_i(\sigma)=m \right) \\
\leq &\ {2 \over n \delta^2}  \sum_{m:\ m \leq n/2} {1 \over m} + {2 \over n} \sum_{m: \ m>n/2} \PP\left(\sigma:\ l_i(\sigma)=m\right) \\ \leq &\ {2 \ln n \over n \delta^2} + {2 \over n}.\endsplit$$
Therefore,
$$\EE c=\sum_{i=1}^n \EE \left(l_i^{-1}\right) \ \leq \ 2+ {2 \ln n \over \delta^2},$$
as desired.
{\hfill \hfill \hfill} \qed
\enddemo

\subhead (2.3) Proof of Theorem 1.3 \endsubhead
Using Lemma 2.2 and the Markov inequality, we obtain
$$\PP\left(\sigma: \ c(\sigma) > 4 + 4 \delta^{-2} \ln n \right) \ \leq \ {\EE c \over 4 + 4 \delta^{-2} \ln n} \ \leq \ {1 \over 2}.$$
Therefore,
$$\sum \Sb \sigma \in S_n: \\ c(\sigma) \leq 4 + 4 \delta^{-2} \ln n \endSb \prod_{i=1}^n a_{i \sigma(i)} =
\left(\per A\right) \PP\left(\sigma: \ c(\sigma) \leq 4 + 4 \delta^{-2} \ln n \right) \ \geq \ {1 \over 2} \per A,$$
as desired.
{\hfill \hfill \hfill} \qed

\subhead (2.4) Proof of Theorem 1.2 \endsubhead
Clearly, 
$$\ham A \ \leq \ \per A.$$
Let
$$\Sigma=\Bigl\{ \sigma \in S_n:\ c(\sigma) \leq 4 + 4 \delta^{-2} \ln n \Bigr\},$$
so by Theorem 1.3
$$\sum_{\sigma \in \Sigma} \prod_{i=1}^n a_{i \sigma(i)} \ \geq \ {1 \over 2} \per A. \tag2.4.1$$
We construct a map $\phi: \Sigma \longrightarrow H_n$, where $H_n \subset S_n$ is the set of all Hamiltonian cycles in $S_n$,
as follows. For a permutation $\sigma \in \Sigma$, we pick the largest element of each cycle and order the cycles in the increasing order of those elements. Let $j_1 < j_2 < \ldots < j_r$ be those elements. Then we patch the cycles into a Hamiltonian cycle $\tau$: we replace 
$$i_1 \rightarrow j_1,\quad i_2 \rightarrow j_2, \quad \ldots, \quad i_r \rightarrow j_r$$
by 
$$i_1 \rightarrow j_2, \quad i_2 \rightarrow j_3, \quad \ldots, \quad i_r \rightarrow j_1.$$
Because of (1.1.1) we have 
$$\PP(\tau) \ \geq \ \delta^{c(\sigma)} \PP(\sigma) \ \geq \ \delta^{4 + 4 \delta^{-2} \ln n} \PP(\sigma).$$
On the other hand, any choice of $j_1 < j_2 < \ldots < j_r$ in 
$$\tau: \ j_1 \rightarrow \cdots \rightarrow j_2 \rightarrow \ldots \rightarrow j_r \rightarrow j_1$$
recovers at most one permutation $\sigma \in \phi^{-1}(\tau)$.
Since every cycle $\tau \in H_n$ corresponds to at most 
$$\sum_{r=0}^{4 + 4\delta^{-2} \ln n} {n \choose r}$$
permutations $\sigma \in \Sigma$, the proof follows.
{\hfill \hfill \hfill} \qed

\head 3. Proof of Theorem 1.8 \endhead 

\proclaim{(3.1) Lemma} Let $\Delta_{m,n}$ be the probability space as in Theorem 1.8. Let us define an $n$-variate polynomial of degree $m$ by 
$$\split f\left(x_1, \ldots, x_n\right) =&\sum \Sb a \in \Delta_{m,n} \\ a=\left(\alpha_1, \ldots, \alpha_n\right) \endSb 
\PP(a) x_1^{\alpha_1} \cdots x_n^{\alpha_n} \\
=&\sum \Sb a \in \Delta_{m,n} \\ a=\left(\alpha_1, \ldots, \alpha_n\right) \endSb 
{w\left(\alpha_1, \ldots, \alpha_n\right) \over \alpha_1! \cdots \alpha_n!} x_1^{\alpha_1} \cdots x_n^{\alpha_n}.
\endsplit$$
Let us define
$$f_i={\partial f \over \partial x_i} \quad \text{for} \quad i=1, \ldots, n.$$
Then, for all $1 \leq i, j \leq n$ and all $x_1, \ldots, x_n \geq 0$ we have 
$$f_i\left(x_1, \ldots, x_n \right) \ \leq \ \delta^{-1} f_j\left(x_1, \ldots, x_n \right).$$
\endproclaim
\demo{Proof}
For $\left(\alpha_1, \ldots, \alpha_n \right) \in \Delta_{m-1, n}$, let us define
$$u_i\left(\alpha_1, \ldots, \alpha_n \right) =w\left(\alpha_1, \ldots, \alpha_{i-1}, \alpha_i +1, \alpha_{i+1}, \ldots, \alpha_n \right).$$
Then
$$\split u_i\left(\alpha_1, \ldots, \alpha_n\right) =&w\left(\alpha_1, \ldots, \alpha_{i-1}, \alpha_i +1, \alpha_{i+1}, \ldots, \alpha_n \right) \\
\ \leq \ &\delta^{-1} w\left(\alpha_1, \ldots, \alpha_{j-1}, \alpha_j +1, \alpha_{j+1}, \ldots, \alpha_n \right)=
\delta^{-1} u_j \left(\alpha_1, \ldots, \alpha_n \right). \endsplit $$
We have
$$\split f_i\left(x_1, \ldots, x_n\right) =&\sum \Sb \left(\alpha_1, \ldots, \alpha_n\right) \in \Delta_{m,n} \\ \alpha_i >0 \endSb
{w\left(\alpha_1, \ldots, \alpha_n\right) \over \alpha_1! \cdots \alpha_{i-1}! \left(\alpha_i-1\right)! \alpha_{i+1}! \cdots \alpha_n!}
\\ &\quad \times
x_1^{\alpha_1} \cdots x_{i-1}^{\alpha_{i-1}} x_i^{\alpha_i -1} x_{i+1}^{\alpha_{i+1}} \cdots x_n^{\alpha_n} 
\\ =& \sum \Sb \left(\alpha_1, \ldots, \alpha_n\right) \in \Delta_{m-1,n}\endSb {u_i\left(\alpha_1, \ldots, \alpha_n\right) 
\over \alpha_1! \cdots \alpha_n!} x_1^{\alpha_1} \cdots x_n^{\alpha_n} \\ 
\leq \ &  \delta^{-1} \sum \Sb \left(\alpha_1, \ldots, \alpha_n\right) \in \Delta_{m-1,n}\endSb {u_j\left(\alpha_1, \ldots, \alpha_n\right) 
\over \alpha_1! \cdots \alpha_n!} x_1^{\alpha_1} \cdots x_n^{\alpha_n} \\
=&\ \delta^{-1} f_j\left(x_1, \ldots, x_n \right).\endsplit$$
{\hfill \hfill \hfill} \qed
\enddemo

Since the polynomial $f$ is homogeneous of degree $m$, by Euler's formula we have 
$$f\left(x_1, \ldots, x_n\right)={1 \over m} \sum_{i=1}^n x_i f_i \left(x_1, \ldots, x_n\right).$$

\proclaim{(3.2) Lemma} Let $\Delta_{m,n}$ be the probability space as in Theorem 1.8 and let 
$f$ be the polynomial as defined in Lemma 3.1.
Then
$$f\left(e^t, 1, \ldots, 1\right) \ \leq \ \left({e^t + (n-1) \delta \over 1+(n-1)\delta} \right)^m \quad \text{for all} \quad t\geq 0.$$
\endproclaim
\demo{Proof}
We have 
$${d \over dt} f\left(e^t, 1, \ldots, 1\right)=e^t f_1\left(e^t, 1, \ldots, 1\right).\tag3.2.1$$
Using Euler's formula and Lemma 3.1, we obtain
$$\split f\left(e^t, 1, \ldots, 1\right) =&{e^t \over m} f_1\left(e^t, 1, \ldots, 1\right) + {1 \over m} \sum_{i=2}^n f_i\left(e^t, 1, \ldots, 1\right)\\ \geq\ &{e^t \over m} f_1\left(e^t, 1, \ldots, 1\right) +{n-1 \over m} \delta f_1\left(e^t, 1, \ldots, 1\right) \\
=&{e^t + (n-1)\delta \over m} f_1\left(e^t, 1, \ldots, 1 \right) \endsplit$$
and hence
$$f_1\left(e^t, 1, \ldots, 1\right) \ \leq \ {m \over e^t +(n-1)\delta} f\left(e^t, 1, \ldots, 1\right). \tag3.2.2$$
Let us denote
$$F(t)=f\left(e^t, 1, \ldots, 1\right).$$
Then $F(0)=1$ and combining (3.2.1)--(3.2.2), we obtain
$${d \over dt } F(t) \ \leq \ {m e^t \over e^t +(n-1) \delta} F(t).$$
Hence
$${d \over dt} \ln F(t)\  \leq \ {m e^t \over e^t +(n-1) \delta}$$ and
$$\ln F(t) \ \leq \ \ln F(0) + \int_0^t {m e^s \over e^s +(n-1) \delta} \ ds = m \ln {e^t +(n-1) \delta \over 1+(n-1) \delta}.$$
The proof now follows.
{\hfill \hfill \hfill} \qed
\enddemo

\subhead (3.3) Proof of Theorem 1.8 \endsubhead Let us define the polynomial $f$ as in Lemma 3.1. Considering the coordinate $\alpha_1$ of $a=\left(\alpha_1, \ldots, \alpha_n \right)$ as a random variable on $\Delta_{m,n}$ and using Lemma 3.2, we obtain 
$$\EE e^{t \alpha_1} =f\left(e^t, 1, \ldots, 1\right) \ \leq \ \left({e^t +(n-1) \delta \over 1+(n-1)\delta} \right)^m.$$
In particular, for $t=\ln \ln n$, we obtain
$$\EE e^{t \alpha_1} \ \leq \ \left( {\ln n +(n-1) \delta \over 1+(n-1)\delta} \right)^m \ \leq \ 
\left(1+ {\ln n \over n \delta}\right)^m \ \leq \ \exp\left\{ m \ln n\over n \delta \right\}.$$
By the Markov inequality,
$$\split \PP\left( \alpha_1 \ > \ {3m \ln n  \over \delta n \ln \ln n} \right) = &\PP\left( e^{t \alpha_1} \ > \ \exp\left\{ {3 m \ln n \over \delta n} \right\}\right) \\ \leq \ &\exp\left\{-{3m \ln n \over \delta n}\right\} \EE e^{t \alpha_1} 
\ \leq \ \exp\left\{-{2m \ln n \over \delta n}\right\} \\
\leq \ &{1 \over n^2}. \endsplit$$
Similarly,
$$\PP\left(\alpha_i \ > \ {3m \ln n \over \delta n \ln \ln n} \right) \ \leq \ {1 \over n^2} \quad \text{for} \quad i=1, \ldots, n$$
and the proof follows.
{\hfill \hfill \hfill} \qed

\head 4. Proofs of Theorems 1.5 and 1.7 \endhead

\subhead (4.1) Proof of Theorem 1.5 \endsubhead We use Theorem 1.8. Let $\Delta_{n, n}$ be the set of all 
non-negative integer $n$-vectors $a=\left(\alpha_1, \ldots, \alpha_n\right)$ such that $\alpha_1 +\ldots +\alpha_n=n$.
We introduce a probability measure on $\Delta_{n, n}$ as follows. For $a \in \Delta_{n,n}$, 
$a=\left(\alpha_1, \ldots, \alpha_n\right)$, let $\Pi(a)$ be the set of closed walks $\pi$ of length $n$ such that 
$$\deg_i(\pi)=\alpha_i \quad \text{for} \quad i=1, \ldots, n.$$
We let
$$\PP(a) =\left(\tr A^n \right)^{-1} \sum_{\pi \in \Pi(a)} \wei(\pi).$$
Let 
$$a=\left(\alpha_1, \ldots, \alpha_n \right)$$ and suppose that $\alpha_1 > 0$. 
We let 
$$b=\left(\alpha_1-1, \alpha_2+1, \alpha_3, \ldots, \alpha_n \right)$$ and compare $\PP(a)$ and $\PP(b)$. 
For each closed walk $\pi \in \Pi(a)$ we construct $\alpha_1$ closed walks $\rho_i$ as follows. Let 
$$\pi: \quad i_1 \rightarrow i_2 \rightarrow \ldots \rightarrow i_n \rightarrow i_1.$$
For each of the $\alpha_1$ occurrences of $i_k=1$ we 
$$\text{replace} \quad i_{k-1} \rightarrow 1 \rightarrow i_{k+1} \quad \text{by} \quad i_{k-1} \rightarrow 2 \rightarrow i_{k+1}$$
(with the obvious adjustment if $k=1$). For every closed walk $\rho_i$ so obtained, we have
$$\wei\left(\rho_i\right) \ \leq \ \delta^{-2} \wei(\pi).$$
Moreover, each closed walk $\rho \in \Pi(b)$ can be obtained in this way from precisely $\alpha_2+1$ closed walks $\pi$ (we apply the reverse operation to $\rho$ in $\alpha_2+1$ positions).
Hence 
$$\left(\alpha_2 +1\right) \PP(b) \ \leq \ \delta^{-2} \alpha_1 \PP(a).$$
Denoting
$$w\left(a\right) =\alpha_1! \cdots \alpha_n! \PP(a) \quad \text{where} \quad 
a=\left(\alpha_1, \ldots, \alpha_n \right),$$
we obtain
$$w\left(\alpha_1, \ldots, \alpha_n\right) \ \leq \ \delta^{-2} w\left(\beta_1, \ldots, \beta_n\right) \quad \text{whenever} \quad
\sum_{i=1}^n \left|\alpha_i-\beta_i\right|=2.$$
We use Theorem 1.8.
{\hfill\hfill\hfill} \qed

\remark{(4.2) Remark} Jeff Kahn \cite{Ka13} communicated to the author an alternative, purely combinatorial proof of Theorem 1.5.
\endremark

\subhead (4.3) Proof of Theorem 1.7 \endsubhead We use Theorem 1.8. Let $\Delta_{n-2, n}$ be the set of all non-negative 
integer $n$-vectors $a=\left(\alpha_1, \ldots, \alpha_n\right)$ such that $\alpha_1 + \ldots + \alpha_n =n-2$.
For $a \in \Delta_{n-2, n}$, $a=\left(\alpha_1, \ldots, \alpha_n\right)$, let $T(a)$ be the set of all spanning trees $\tau$ such that 
$\deg_i \tau =\alpha_i +1$ for $i=1, \ldots, n$. We define a probability measure on $\Delta_{n-2,n}$ by
$$\PP(a) =\left(\spt A \right)^{-1} \sum_{\tau \in T(a)} \wei(\tau).$$
Let 
$$a=\left(\alpha_1, \ldots, \alpha_n \right)$$ and suppose that $\alpha_1 > 0$. 
We let 
$$b=\left(\alpha_1-1, \alpha_2+1, \alpha_3, \ldots, \alpha_n \right)$$ and compare $\PP(a)$ and $\PP(b)$. 
For each tree $\tau \in T(a)$ we construct $\alpha_1$ trees $\eta_i \in T(b)$ as follows. There is a unique path $\gamma$ in $\tau$ connecting the vertices $1$ and $2$ and hence there is a unique edge in $\gamma$ adjacent to 1. Therefore, there is a set $S$ of 
exactly $\alpha_1$ vertices $i$ such that $\{1, i\}$ is an edge of $\tau$ and $i \notin \gamma$. Furthermore, for every $i \in S$ the vertices $i$ and $2$ are not connected by an edge in $\tau$, as that would have resulted in a cycle in $\tau$. We pick a vertex $i \in S$, remove the edge $\{1, i\}$ from $\tau$ and add the edge $\{2, i\}$ to $\tau$. We get a graph $\eta_i$ with $n-1$ edges which is still connected, because the vertices $1$ and $i$ remain connected via the path $\gamma$ from $1$ to $2$ and then by the edge $\{2, i\}$. Hence $\eta_i$ is a spanning tree. We have
$$\wei\left(\eta_i\right) \ \leq \ \delta^{-1} \wei(\tau).$$
Moreover, each tree $\eta \in T(b)$ is obtained from precisely $\alpha_2+1$ trees this way, as we can apply the reverse procedure to $\eta$ in $\alpha_2+1$ ways.
Hence
$$\left(\alpha_2+1\right) \PP(b) \ \leq \ \delta^{-1} \alpha_1 \PP(a)$$
and the proof is finished as in Section 4.1.
{\hfill \hfill \hfill} \qed

\remark{(4.4) Remark} Jeff Kahn \cite{Ka13} suggested that an alternative proof of Theorem 1.7 can be based on the Aldous - Broder algorithm for generating a random spanning tree, see \cite{Al90}, \cite{B789} and \cite{MS99} for the weighted version.
\endremark

\head Acknowledgment \endhead

The author is indebted to Jeff Kahn for suggesting an alternative combinatorial proof of Theorems 1.5 and  pointing out to connections of Theorem 1.7 with the Aldous - Broder algorithm.

\enddocument
\end